\documentclass[12pt]{article}
\begin{document}
\title{A NEW INTEGRAL TRANSFORM}
\author{B.G. Sidharth\\
International Institute for Applicable Mathematics \& Information Sciences\\
Hyderabad (India) \& Udine (Italy)\\
B.M. Birla Science Centre, Adarsh Nagar, Hyderabad - 500 063 (India)}
\date{}
\maketitle
\begin{abstract}
Using Bauer's expansion and properties of spherical Bessel and Legender functions, we deduce a new transform and briefly indicate its use.
\end{abstract}
Using properties of spherical Bessel and Legender functions, we would now like to deduce a new Integral Transform.
Our starting point is Bauer's expansion (this and a few other known results quoted can be obtained from refs. \cite{wab}-\cite{bat}): 
$$e^{\imath zt} = \sum^{\infty} (2n+1) \imath^n P_n (t) j_n (z)$$
Using the orthogonalty of the spherical Bessel and Legender functions, viz., the relations
$$\int^{\infty}_{-\infty} j_n (z) j_m (z) dz = \int^{\infty}_{-\infty} J_{n+\frac{1}{2}} (z) J_{m+\frac{1}{2}} (z)\frac{dz}{z} = 0$$
if 
$$m \ne n$$
and
$$= \frac{2}{(2n+1)} \, \mbox{if}\, m = n$$
$$\int^{+1}_{-1} P_n (t) P_m (t) dt = 0$$
if 
$$m \ne n,$$
and
$$= \frac{2}{2n+1} \, \mbox{if} \, m = n$$
in Bauer's expansion we get
\begin{equation}
\int^1_{-1} e^{\imath z t} P_n (t) dt = 2 \imath^n j_n(z)\label{e1}
\end{equation}
\begin{equation}
\int^\infty_{-\infty} j_n (z) e^{\imath z t} dz = 2 \imath^n P_n(t)\label{e2}
\end{equation}
We consider that $z$ is real. 
We would also need the following
$$J_\nu (z) = \frac{(\frac{1}{2}z)^\nu}{\Gamma (\nu + \frac{1}{2})\Gamma (\frac{1}{2})} \int^\pi_0 \cos (z \cos \Theta ) \sin^{2 \nu} \Theta d\Theta$$
$$\left(Re (\nu + \frac{1}{2}) > 0 \right) \quad \nu = n + \frac{1}{2}$$
Whence,
$$J_\nu (-z) = (-1)^\nu J_\nu (z) = (-1)^n \imath J_\nu (z)$$
or 
\begin{equation}
\dot {j}_n (-z) = (-1)^n \dot{j}_n (z)\label{ec}
\end{equation}
Let us consider a function $g (z)$ which can be expanded as an infinite linear combination of spherical Bessel functions, on the lines of Neumann's expansion in terms of ordinary Bessel functions. This can be done because of the orthogonality relations above. Similarly, we will also use the known expansion in terms of Legender functions. Thus we have,
$$g(z) = c_n j_n (z)$$
$$= \sum c_n (2\imath^n)^{-1} \int^1_{-1} e^{\imath zt} P_n (t) dt$$
or,
\begin{equation}
g(z) = \int^1_{-1} f(t) e^{\imath zt} dt\label{e3}
\end{equation}
where
$$f(t) = \sum \bar{c}_n P_n (t), (\bar{c}_n = c_n (2\imath^n)^{-1})$$
$$= \sum \bar{c}_n (2\imath^n)^{-1} \int^\infty_{-\infty} j_n (z) e^{\imath zt} dz$$
or
$$f(t) = \int \frac{1}{4} \sum c_n (-1)^n j_n (z) e^{\imath zt} dz$$
$$= \frac{1}{4} \int^\infty_{-\infty} \sum c_n j_n (-z) e^{-\imath (-z)t} dz$$
wherein we have used (\ref{ec}), or,
\begin{equation}
f(t) = \frac{1}{4} \int^\infty_{-\infty} g(y) e^{-\imath yt} dy\label{e4}
\end{equation}
In deducing (\ref{e3}) and (\ref{e4}), we have used (\ref{e1}) and (\ref{e2}), and the summations are infinite. Moreover we assume that for $f(t)$ and $g(z)$ derivatives of all orders exist over their domains.\\ 
So finally,
\begin{equation}
g (z) = \frac{1}{4} \int^{+1}_{-1} \int^\infty_{-\infty} g(y) e^{\imath (z-y)t} dydt\label{e5}
\end{equation}
The relations (\ref{e3}), (\ref{e4}) and (\ref{e5}) are the desired new relations. As an application, let us consider the differential equation,
\begin{equation}
L_{op} g(z) = h (z),\label{ed}
\end{equation}
where $L_{op}$ is a linear differential operator. Using (\ref{e5}) in (\ref{ed}), we get,
$$L_{op} g (z) = F(\frac{d}{dz}) g(z) = \frac{1}{4} \int^{+1}_{-1} \int^{\infty}_{-\infty} F(\imath t) g(y) e^{\imath (z-y)t} dy dt = h (z)$$
or
$$L_{op} g(z) = A \int^{+1}_{-1} f(t) F(\imath t)e^{\imath zt} dt = h (z)$$
\begin{equation}
= \int^{+1}_{-1} \hat {h} (t) e^{\imath zt} dt\label{e6}
\end{equation}
where we have used (\ref{e3}),
$$h (z) = \sum d_n j_n (z), \quad \hat {h} (t) = \sum \bar{d}_n P_n (t)$$
$$\bar{d}_n = d_n (2 \imath^n)^{-1}$$
So we get
$$f(t) F(\imath t) = \hat {h} (t)$$
As $\hat {h} (t) \, \mbox{and} \, F(\imath t)$ are known so is $f(t)$ known and therefore also $g(z)$. Infact
$$f(t) = \frac{1}{4} \int^{\infty}_{-\infty} g(y) e^{-\imath yt} dy = \sum \bar{c}_n P_n (\imath t)$$
so that
$$g (z) = \sum c_n j_n (z), \bar{c}_n = (2 \imath^n)^{-1} c_n$$
{\large {\bf Remarks:}}\\ \\
1. We note that Neumann's expansion alluded to applies for any analytical function $g(z)$:
$$g(z) = \sum^{\infty}_{n} b_n J_n (z)$$
However the expansion in (\ref{e3}) is in terms of Spherical Bessel functions. As mentioned such an expansion can always be justified, as in the case of the Legender polynominal expansion of any function $f(t)$ given in (\ref{e4}), by using the orthogonalty properties of the $j_n(z)$ and $P_n (t)$ given above.\\
2. The above consideration in relation (\ref{e5}) is to be distinguished from the so called Hankel transform. Further, it must be noted that the domains of integration in (\ref{e3}), (\ref{e4}) and (\ref{e5}) are $(-1,1)$ for $t$ and $(-\infty , \infty)$ for $z$.
 

\begin{thebibliography}{99}
\bibitem {wab} G.N. Watson, ``Theory of Bessel Functions'', Cambridge University Press, Cambridge, 1958.
\bibitem {whit} E.T. Whittaker and G.N. Watson, ``A Course of Modern Analysis'', Cambridge University Press, Cambridge, 1962.
\bibitem {mase} P.M. Morse and H. Feshbach, ``Methods of Theoretical Physics'', Vol.2, McGraw Hill, New York, 1958.
\bibitem {cop} E.T. Copson, ``Theory of Functions of a Complex Variable'', University Press, London, 1935.
\bibitem {bat} H. Bateman, ``Higher Transcendental Functions'', Vol2. McGraw Hill, New York, 1953.
\end{thebibliography}
\end{document}